\documentclass[10pt]{article}
\usepackage[nofoot,twoside,paperheight=235mm,paperwidth=160mm]{geometry}
\usepackage{amsmath,amssymb,amsthm,amsfonts}\usepackage[hidelinks]{hyperref}
\usepackage[english]{babel}\usepackage{graphicx}
\theoremstyle{plain}\newtheorem{theorem}{Theorem}\newtheorem{proposition}{Proposition}
\theoremstyle{definition}
\theoremstyle{remark}\newtheorem{remark}{Remark}
\usepackage{fancyhdr}
\long\def\symbolfootnote[#1]#2{\begingroup\def\thefootnote{\fnsymbol{footnote}}\footnote[#1]{#2}\endgroup}
\def\title#1{{\Large\bf  \begin{center} #1 \vspace{0pt} \end{center}  } \smallskip}
\def\authors#1{{\bf \begin{center} #1 \vspace{0pt} \end{center} } \smallskip}
\def\institution#1{{\sl \begin{center} #1 \vspace{0pt} \end{center} } }
\def\keywords#1{\bigskip \par\noindent{\bf Keywords: }#1\par}
\def\AMS#1{\par\noindent{\bf AMS subject classifications: }#1\par}
\newcommand{\acknowledgements}{\par{\noindent\bf Acknowledgements: }}
\begin{document}\fancyhead[LE,RO]{{\small\arabic{page}}}
\fancyhead[LO]{{\small Proceedings of the 19th EYSM in Prague 2015}}
\fancyfoot{}\thispagestyle{empty}
\renewenvironment{abstract}{\bigskip\noindent\begin{minipage}{\textwidth}\setlength{\parindent}{15pt}\paragraph{Abstract:}}{\end{minipage}}

\fancyhead[RE]{{\small Rossi -- High Energy Random Eigenfunctions on $\mathbb S^d$}}

\title{On the High Energy Behavior of Nonlinear Functionals of Random Eigenfunctions on $\mathbb S^d$}

\authors{Maurizia Rossi\symbolfootnote[1]{Corresponding author:
rossim@mat.uniroma2.it}}

\institution{Dipartimento di Matematica, Universit\`a di Roma ``Tor Vergata''}


\begin{abstract}
In this short survey we recollect some of the recent results
on the high energy behavior (i.e., for diverging sequences of eigenvalues)
of nonlinear functionals of Gaussian eigenfunctions on the $d$-dimensional
sphere $\mathbb S^d$, $d\ge 2$.
We present a quantitative Central Limit Theorem for a class of functionals
whose Hermite rank is two, which includes in particular
 the empirical measure of excursion sets in the non-nodal case. Concerning the   nodal case,
 we recall a CLT result for the defect on $\mathbb S^2$.
The key tools are
both, the asymptotic analysis of moments of all order for Gegenbauer polynomials,
and  so-called Fourth-Moment theorems.
\end{abstract}

\keywords{Gaussian eigenfunctions, High Energy Asymptotics, Quantitative Central
Limit Theorems, Excursion Volume, Gegenbauer polynomials.  }

\AMS{60G60, 42C10, 60D05, 60B10.}

\section{Introduction}

Let us consider a compact Riemannian manifold $(\mathcal M, g)$
and denote by $\Delta_{\mathcal M}$ its Laplace-Beltrami operator. There
exists a sequence of eigenfunctions $\lbrace f_j \rbrace_{j\in \mathbb N}$
and a corresponding non-decreasing sequence of eigenvalues $\lbrace E_j \rbrace_{j\in \mathbb N}$
\begin{equation*}
\Delta_{\mathcal M}\,  f_j + E_j f_j = 0\ ,
\end{equation*}
such that $\lbrace f_j \rbrace_{j\in \mathbb N}$ is a complete orthonormal basis of $L^2(\mathcal M)$,
the space of
square integrable measurable functions on $\mathcal M$. One is interested in the
high energy behavior i.e., as $j\to +\infty$, of eigenfunctions $f_j$, related to the
geometry of both  \emph{level sets} $f_j^{-1}(z)$ for $z\in \mathbb R$,
and connected components of their complement $\mathcal M \setminus f_j^{-1}(z)$.
One can investigate e.g. the Riemannian volume of these domains: a
 quantity that can be formally written  as a \emph{nonlinear functional} of $f_j$.

The nodal case corresponding to $z=0$ has received  great attention
 (for  motivating details see \cite{wigsurvey}).

At least for ``generic'' chaotic \emph{surfaces} $\mathcal M$, Berry's Random Wave Model allows
to compare the   eigenfunction $f_j$ to a ``typical''
instance of an isotropic, monochromatic \emph{random} wave with wavenumber
$\sqrt{E_j}$ (see \cite{wigsurvey}).
In view of this, much effort has been first devoted  to $2$-dimensional manifolds such as the  torus $\mathbb T^2$ (see e.g. \cite{AmP}) and the sphere $\mathbb S^2$ (see e.g.  \cite{vale2}, \cite{vale3},  \cite{Nonlin},
 \cite{Wig}). Spherical random fields have  attracted  a growing interest, as  they  model several data sets
in Astrophysics and Cosmology, e.g. on Cosmic Microwave Background (\cite{giodom}).

More recently  random eigenfunctions on higher dimensional manifolds have been investigated: e.g.  on the hyperspheres  (\cite{maudom}).
\subsection{Random eigenfunctions on $\mathbb S^d$}\label{subsec}

Let us fix some probability space $(\Omega, \mathcal F, \mathbb P)$, denote by $\mathbb E$ the corresponding expectation and by $\mathbb S^d\subset \mathbb R^{d+1}$  the unit $d$-dimensional sphere ($d\ge 2$);
$\mu_d$ stands for the Lebeasgue measure of the hyperspherical surface.
By  real random field  on $\mathbb S^d$ we mean a real-valued measurable map defined on $(\Omega\times \mathbb S^d, \mathcal F\otimes \mathcal B(\mathbb S^d))$, where $\mathcal B(\mathbb S^d)$  denotes the Borel $\sigma$-field on $\mathbb S^d$.
Recall that the eigenvalues of the Laplace-Beltrami operator $\Delta_{\mathbb S^d}$ on $\mathbb S^d$ are
 integers of the form $ -\ell(\ell+d-1)=: - E_\ell$, $\ell \in \mathbb N $.

The $\ell$-th random eigenfunction $T_\ell$ on $\mathbb S^d$
is the (unique) centered, isotropic real Gaussian  field  on $\mathbb S^d$ with
covariance function
\begin{equation*}
K_\ell(x,y):= G_{\ell;d}(\cos \tau(x,y))\quad x,y\in \mathbb S^d\ ,
\end{equation*}
where $G_{\ell;d}$ stands for the $\ell$-th Gegenbauer polynomial normalized in such a way that
$G_{\ell;d}(1)=1$ and $\tau$ is the usual
geodesic distance.  More precisely, setting $\alpha_{\ell;d} := {\ell + \frac{d}2 -1 \choose \ell}$, we have
$
G_{\ell;d} =\alpha_{\ell;d}^{-1}\,P_\ell^{(\frac{d}2-1,\frac{d}2-1)}\ ,
$ where
 $P_\ell^{(\alpha, \beta)}$ denote standard Jacobi polynomials.
By isotropy (see e.g. \cite{giodom}) we mean that for every $g\in SO(d+1)$,
the random fields $T_\ell = (T_\ell(x))_{x\in \mathbb S^d}$ and
$T^g_\ell := (T_\ell(gx))_{x\in \mathbb S^d}$ have the same law in the sense
of finite-dimensional distributions. Here $SO(d+1)$ denotes the group
of real $(d+1)\times (d+1)$-matrices $A$ such that $AA'=I$ the identity matrix and $\text{det} A =1$.

Random eigenfunctions naturally arise as they are the Fourier components of those isotropic random fields on $\mathbb S^d$ whose sample paths belong to  $L^2(\mathbb S^d)$.

Let us consider now
 functionals of  $T_\ell$
of the form
\begin{equation}\label{Mm}
S_\ell(M) := \int_{\mathbb S^d} M(T_\ell(x))\,dx\ ,
\end{equation}
where $M:\mathbb R \to \mathbb R$ is some measurable function such that $\mathbb E[M(Z)^2] < +\infty$,
$Z\sim \mathcal N(0,1)$ a standard Gaussian r.v.
In particular, if $M(\cdot) = 1(\cdot > z)$ is the indicator function of the interval $(z, +\infty)$ for $z\in \mathbb R$,
then (\ref{Mm}) coincides with the  empirical measure $S_\ell(z)$ of the $z$-\emph{excursion set}
$
A_\ell(z) := \lbrace x\in \mathbb S^d : T_\ell(x) > z \rbrace
$.
\subsection{Aim of the survey}
We first present a quantitative CLT as $\ell\to +\infty$ for nonlinear functionals $S_\ell(M)$ in (\ref{Mm}) on $\mathbb S^d$, $d\ge 2$,  under the assumption that $\mathbb E[M(Z) H_2(Z)]\ne 0$, where $H_2(t) := t^2 -1$  is the second Hermite polynomial.

 For instance the above condition is fullfilled by the empirical measure  $S_\ell(z)$ of $z$-excursion sets for $z\ne 0$. For the nodal case which corresponds to the defect
\begin{equation}\label{defect}
D_\ell := \int_{\mathbb S^d} 1(T_\ell(x) > 0)\,dx -  \int_{\mathbb S^d} 1(T_\ell(x) < 0)\,dx\ ,
\end{equation}
 we present a CLT for $d=2$. Quantitative CLTs for  $D_\ell$ on $\mathbb S^d$, $d \ge 2$, will be treated in a forthcoming paper.

We refer to \cite{Def}, \cite{Nonlin} and \cite{maudom} for the spherical case $d=2$ and to \cite{maudom} for all higher dimensions.
The mentioned results rely  on both, the asymptotic analysis of moments of all order
for Gegenbauer polynomials, and Fourth-Moment theorems (see \cite{noupe}, \cite{simon}).

\section{High energy behavior  via chaos expansions}

For a function $M:\mathbb R \to \mathbb R$ as in (\ref{Mm}), the r.v. $S_\ell(M)$ admits the chaotic expansion
\begin{equation}\label{series}
S_\ell(M) = \sum_{q=0}^{+\infty} \frac{J_q(M)}{q!} \int_{\mathbb S^d} H_q(T_\ell(x))\,dx
\end{equation}
 (see \cite{noupe}) in $L^2(\mathbb P)$ (the space of finite-variance r.v.'s),
where $H_q$ is the $q$-th Hermite polynomial (see e.g. \cite{szego}) and $J_q(M) := \mathbb E[M(Z)H_q(Z)]$,
 $Z\sim \mathcal N(0,1)$.  We have  $\mathbb E[S_\ell(M)] = J_0(M) \mu_d$;   w.l.o.g. $J_0(M)=0$.

The main idea is first to investigate the asymptotic behavior of each chaotic projection,
i.e. of each (centered) r.v. of the form
\begin{equation}\label{hell}
h_{\ell;q,d}:=\int_{\mathbb S^d} H_q(T_\ell(x))\,dx
\end{equation}
and then deduce the asymptotic behavior of the whole series (\ref{series}).
Note that $h_{\ell;1,d} = 0$, as $T_\ell$ has zero mean on $\mathbb S^d$.
By the symmetry property of Gegenbauer polynomials (\cite{szego}),
from now on we can restrict ourselves to even multiples
$\ell$, for which some straightforward computations yield
\begin{equation}\label{var}
\text{Var}[h_{\ell;q,d}]= 2q! \mu_d \mu_{d-1} \int_0^{\pi/2}
G_{\ell;d}(\cos \vartheta)^q (\sin \vartheta)^{d-1}\, d\vartheta\ .
\end{equation}
\subsection{Asymptotics for moments of Gegenbauer polynomials}
The proof of the following is in \cite{Def}, \cite{Nonlin} for $d=2$ and in \cite{maudom} for $d\ge 3$.
\begin{proposition}
\label{varianza} As $\ell \rightarrow \infty ,$ for $d=2$ and $q=3$ or $q\ge 5$ and for $d,q\ge 3$,
\begin{equation}  \label{int1}
\int_{0}^{\frac{\pi }{2}}G_{\ell ;d}(\cos \vartheta )^{q}(\sin \vartheta
)^{d-1}\,d\vartheta =\frac{c_{q;d}}{\ell ^{d}}(1+o(1))\ .
\end{equation}
The constants $c_{q;d}$ are given by the formula
\begin{equation}
c_{q;d}:=\left(2^{\frac{d}{2} - 1}\left (\frac{d}2-1 \right)!\right)^q\int_{0}^{+\infty
}J_{\frac{d}{2}-1}(\psi )^{q}\psi ^{-q\left( {\textstyle\frac{d}{2}}%
-1\right) +d-1}d\psi\ ,  \label{cq}
\end{equation}%
where $J_{\frac{d}{2}-1}$ is the Bessel function (\cite{szego})
of order $\frac{d}{2}-1$.
The r.h.s. integral in (\ref{cq}) is absolutely convergent for any
pair $(d,q)\neq (2,3), (3,3)$ and conditionally convergent for $d=2, q=3$ and $d=q=3$.
Moreover for $c_{4;2}:=\frac{3}{2\pi^2}$
\begin{equation}
\int_{0}^{\frac{\pi }{2}}G_{\ell ;2}(\cos \vartheta )^{4}\sin \vartheta
\,d\vartheta = \ c_{4;2}\frac{\log \ell }{\ell ^{2}}(1 + o(1))\
.  \label{q=4d=2}
\end{equation}%
\end{proposition}
From \cite{szego}, as $\ell\to +\infty$,
\begin{equation}
\int_{0}^{\frac{\pi }{2}}G_{\ell ;d}(\cos \vartheta )^{2}(\sin \vartheta
)^{d-1}\,d\vartheta = 4\mu_d
 \mu_{d-1}\frac{c_{2;d}}{\ell ^{d-1}}(1+o(1))\ ,\  c_{2;d} :=
\frac{(d-1)!\mu _{d}}{4 \mu_{d-1}}\ . \label{q=2}
\end{equation}%
Clearly for any $d,q\geq 2$,  $c_{q;d}\ge 0$ and $c_{q;d}>0$ for all even $q$. Moreover we can give explicit expressions for $c_{3;2}, c_{4;2}$ and $c_{2;d}$ for any $d\ge 2$.
We conjecture that the above strict inequality holds for every pair $(d,q)$,
and leave this issue as an open question for future research.
\subsection{Fourth-Moment Theorems for chaotic projections}
Let us recall the
usual Kolmogorov $d_{K}$, total variation $d_{TV}$ and Wasserstein $d_{W}$
distances between r.v.'s $X,Y$: for $\mathcal D \in \lbrace K, TV, W \rbrace$
$$
d_{\mathcal D}(X,Y) :=\sup_{h\in  H_{\mathcal D}}\left\vert {\mathbb{E}}[h(X)]-{%
\mathbb{E}}[h(Y)]\right\vert \text{ ,}  \label{prob distance}
$$
where $H_{K} = \lbrace  1(\cdot \le z), z\in \mathbb R \rbrace$,
 $H_{TV} = \lbrace  1_A(\cdot), A\in \mathcal B(\mathbb R) \rbrace$  and $H_{W}$ is
 the set of Lipschitz functions with Lipschitz
constant one.

The r.v. $h_{\ell;q,d}$ in (\ref{hell}) belongs to the so-called $q$th Wiener chaos. The Fourth-Moment Theorem (\cite{noupe}) states that if $Z\sim \mathcal N(0,1)$, for $\mathcal D \in \lbrace K, TV, W \rbrace$ we have
\begin{equation}\label{4th}
d_{\mathcal D} \left( \frac{h_{\ell ;q,d}}{\sqrt{\text{Var}[h_{\ell ;q,d}]}},Z\right) \le C_{\mathcal D} (q) \sqrt{\frac{\text{cum}_4(h_{\ell ;q,d})}{\text{Var}[h_{\ell ;q,d}]^2}}\ ,
\end{equation}
where $C_{\mathcal D}(q)>0$ is some explicit constant and  $\text{cum}_4(h_{\ell ;q,d})$ is the fourth cumulant of the r.v. $h_{\ell ;q,d}$.
An application of (\ref{4th}) together with upper bounds for cumulants leads to the following result  (see \cite{maudom}).
\begin{theorem}
\label{teo1} For all $d,q\ge 2$ and
${\mathcal{D}}\in \lbrace K, TV, W \rbrace$ we have, as $\ell\to +\infty$,
\begin{equation}\label{bounds}
d_{\mathcal{D}}\left( \frac{h_{\ell ;q,d}}{\sqrt{\mathrm{Var}[h_{\ell ;q,d}]}},%
Z\right) =O \left (\ell^{-\delta(q;d)} (\log \ell)^{-\eta(q;d)}\right )\ ,
\end{equation}%
where $\delta(q;d)\in \mathbb Q$, $\eta(q;d)\in \lbrace -1, 0, 1 \rbrace$ and  $\eta(q;d)=0$ but for $d=2$ and $q=4,5,6$.
\end{theorem}
The exponents  $\delta(q;d)$ and $\eta(q;d)$ can be given explicitly (see \cite{maudom}), turning out  in particular
%
%
that
if $(d,q)\neq (3,3), (3,4),(4,3),(5,3)$ and $%
c_{q;d}>0$,
\begin{equation}\label{convergenza hq}
\frac{h_{\ell ;q,d}}{\sqrt{\text{Var}[h_{\ell ;q,d}]}}\ \mathop{\to}^{\mathcal{L}}\
Z\ ,\qquad \text{as}\ \ell \rightarrow +\infty \ ,
\end{equation}
where from now on, $\mathop{\to}^{\mathcal{L}}$ denotes convergence in distribution and $Z\sim \mathcal N(0,1)$.
\begin{remark}
For $d=2$, the CLT  (\ref{convergenza hq}) was already proved in
\cite{Nonlin}; nevertheless Theorem \ref{teo1} improves the existing bounds
on the rate of convergence to the asymptotic Gaussian distribution.
\end{remark}
\subsection{Quantitative CLTs for Hermite rank $2$ functionals}
Proposition \ref{varianza} states that whenever $M$ is such that  $J_2(M)\ne 0$ in (\ref{series}), i.e.  the functional $S_\ell(M)$ in (\ref{Mm}) has Hermite rank two, then
\begin{equation}\label{=var}
\lim_{\ell\to +\infty} \frac{\text{Var}[S_\ell(M)]}{\text{Var}\left [\frac{J_2(M)}{2}h_{\ell;2,d} \right ]} = 1\ .
\end{equation}
Hence, loosely speaking,  $S_\ell(M)$ and its $2$nd chaotic projection $\frac{J_2(M)}{2}h_{\ell;2,d}$ have the same high energy behaviour.
The main result presented in this survey is the following, whose proof is given in \cite{maudom}.
\begin{theorem}
\label{general} Let $M:\mathbb R\to \mathbb R$ in (\ref{Mm}) be  s.t. $\mathbb{E}\left[
M(Z)H_{2}(Z)\right] =:J_{2}(M)\neq 0$, then
\begin{equation}
d_{W}\left( \frac{S_{\ell }(M)}{\sqrt{%
\mathrm{Var}[S_{\ell }(M)]}},Z\right) =O\left (\ell ^{-\frac{1}{2}} \right ) \
,\qquad \text{as}\ \ell \rightarrow \infty\ ,  \label{sun2}
\end{equation}
where $Z\sim \mathcal N(0,1)$.
In particular, as $\ell\to +\infty$,
\begin{equation}
\frac{S_{\ell }(M)}{\sqrt{\mathrm{Var}[S_{\ell
}(M)]}}\ \mathop{\to}^{\mathcal{L}}\ Z\ .  \label{sun1}
\end{equation}
\end{theorem}
\section{Geometry of  high energy excursion sets}

Consider the empirical measure $S_\ell(z)$ of the $z$-excursion set   $A_\ell(z)$ for $z\in \mathbb R$, as in \S \ref{subsec}.
It is easy to check that in (\ref{series}) ${\mathbb{E}}[S_{\ell }(z)]=\mu _{d}(1-\Phi (z))$ and  for $q\ge 1$, $J_q(1(\cdot >z)) = H_{q-1}(z)\phi(z)$, where
$\Phi$ and $\phi$ denote respectively the cdf and the pdf of the standard Gaussian law. Since $J_2(1(\cdot >z)=z\phi(z)$, Theorem \ref{general} immediately entails that, as $\ell \rightarrow \infty$,
if $z\ne 0$
\begin{equation*}
d_{W}\left( \frac{S_{\ell }(z)-\mu_d(1- \Phi(z))}{\sqrt{\mathrm{Var}[S_{\ell
}(z)]}}, Z\right) =O\left(\ell ^{-\frac{1}{2}}\right)\text{
.}
\end{equation*}
The nodal case $z=0$  requires different arguments: in the chaos expansion for the defect (\ref{defect}) $D_\ell$ only odd chaoses occur but each of them ``contributes'' by Proposition \ref{varianza}.
Asymptotics for the defect variance on $\mathbb S^2$ have been given in \cite{Def}:
\begin{equation*}
\mathrm{Var}[D_\ell] = \frac{C}{\ell^2}( 1 + o(1))\ , \quad \text{as}\  \ell\to +\infty\ ,
\end{equation*}
for $C>\frac{32}{\sqrt{27}}$. Moreover in \cite{Nonlin} a CLT has been proved: as $\ell\to + \infty$,
\begin{equation*}
\frac{D_\ell}{\sqrt{\mathrm{Var}[D_\ell]}}\ \mathop{\to}^{\mathcal{L}}\ Z\ ,  \label{sun1}
\end{equation*}
where $Z\sim \mathcal N(0,1)$. In a forthcoming paper, we will provide quantitative CLTs for the defect on $\mathbb S^d$, $d\ge 2$.

\begin{remark} The volume of excursion sets is just one instance of Lipschtz-Killing curvatures. In the $2$-dimensional
case, these are completed by the  Euler-Poincar\'e
characteristic (\cite{vale2}) and the length of level curves (\cite{AmP},\cite{Wig} for the nodal variances). In forthcoming papers jointly with D.~Marinucci, G.~Peccati and I.~Wigman,
we will investigate the asympotic distribution of the latter on both the sphere $\mathbb S^2$ and the $2$-torus $\mathbb T^2$.
For future research, we would like to  characterize
 the high energy behavior of all Lipschitz-Killing curvatures
on  every  ``nice'' compact manifold.
\end{remark}
\bigskip
\acknowledgements{We thank D.~Marinucci for valuable suggestions, P.~Baldi, S.~Campese and S.~Cipolla for  a careful reading of  an earlier version of this work.

This research is supported by ERC  Grant \emph{Pascal} n.277742.}


\end{document}